\documentclass[submission,copyright,creativecommons]{eptcs}
\providecommand{\event}{ThEdu'23} 
\usepackage{iftex}

\ifpdf
  \usepackage{underscore}         
  \usepackage[T1]{fontenc}        
\else
  \usepackage{breakurl}           
\fi

\usepackage{amsmath}
\usepackage{amssymb,amsthm,amsfonts}
\usepackage{enumerate}
\usepackage{graphicx}

\usepackage{xcolor}
\usepackage{colortbl}
\usepackage{tabularray}
\usepackage[framemethod=TikZ]{mdframed}
\newmdenv[hidealllines=true,backgroundcolor=gray!20,font=\tt]{coqcode}

\usepackage{newunicodechar}
\newunicodechar{ℝ}{$\mathbb{R}$}
\newunicodechar{⇒}{$\Rightarrow$}
\newunicodechar{∧}{$\wedge$}
\newunicodechar{≤}{$\leq$}
\newunicodechar{≥}{$\geq$}
\newunicodechar{ℕ}{$\mathbb{N}$}
\newunicodechar{ℤ}{$\mathbb{Z}$}
\newunicodechar{ℂ}{$\mathbb{C}$}
\newunicodechar{∈}{$\in$}
\newunicodechar{∀}{$\forall$}
\newunicodechar{∞}{$\infty$}
\newunicodechar{∩}{$\cap$}

\usepackage{listings}
\usepackage[frozencache=true,cachedir=minted-cache]{minted} 
\newminted[coqlinenumbers]{coq}{
    linenos=true,
    frame=single,
    framesep=1mm,
    fontsize=\footnotesize,
    numbersep=-10pt
    }
    
\newenvironment{coqlinenumbersstart}[1]
 {\VerbatimEnvironment\begin{coqlinenumbers*}{firstnumber=#1}}
 {\end{coqlinenumbers*}}

\title{Waterproof: Educational Software for Learning\\ How to Write Mathematical Proofs}
\author{Jelle Wemmenhove$^*$ \qquad Dick Arends \qquad Thijs Beurskens \qquad Maitreyee Bhaid \\ Sean McCarren \qquad Jan Moraal \qquad Diego Rivera Garrido \qquad David Tuin \\ Malcolm Vassallo \qquad Pieter Wils \qquad Jim Portegies$^\dag$
\institute{Eindhoven University of Technology\\ Eindhoven, The Netherlands}
\email{$^*$a.j.wemmenhove@tue.nl \quad $^\dag$j.w.portegies@tue.nl} 
}
\def\titlerunning{Waterproof: educational software for learning how to write mathematical proofs}
\def\authorrunning{A.J. Wemmenhove et al.}
\begin{document}
\maketitle

\begin{abstract}
In order to help students learn how to write mathematical proofs, we adapt the \textsc{Coq} proof assistant into an educational tool we call \href{https://impermeable.github.io/}{\emph{Waterproof}}. Like with other interactive theorem provers, students write out their proofs inside the software using a specific syntax, and the software provides feedback on the logical validity of each step.
Waterproof consists of two components:
a custom proof language that allows formal, machine-verified proofs to be written in a style that closely resembles handwritten proofs, and a custom editor that allows these proofs to be combined with formatted text to improve readability. The editor can be used for \textsc{Coq} documents in general, but also offers special features designed for use in education. Student input, for example, can be limited to specific parts of the document to prevent exercises from being accidentally deleted.
Waterproof has been used to supplement teaching the Analysis 1 course at Eindhoven University of Technology (TU/e) for the last four years.
Students started using the specific formulations of proof steps from the custom proof language in their handwritten proofs; the explicit phrasing of these sentences helped to clarify the logical structure of their arguments.
\end{abstract}

\section{Introduction} \label{sec:intro}
Many first-year mathematics undergraduate students struggle with learning how to write mathematical proofs. It is a new skill they often never encountered before, one they nonetheless need to master as one of the key mathematical competencies.
Although students may have difficulty solving the puzzle at the core of many proofs, this is not the main issue \cite{Moore1994,seldenselden2008}.
During various proof stages, students seem to be insufficiently aware what is required of them, and which steps they are, and are not, allowed to perform. 
Students will try to prove $\forall$-statements without introducing a variable; they are quick to assign extra properties to variables obtained from $\exists$-statements; or they swap the order between quantifiers, e.g. in $\varepsilon$-$\delta$ proofs.
Struggling with new concepts is a natural part of the learning process, but too often the confusion about the mechanics underlying mathematical proofs persist beyond the training phase and hinders students' performance in later courses.

The learning process for mathematical proof writing can potentially be improved by the use of proof assistants \cite{french-experience-proof-assistants, Thoma2021}.
Such computer programs, like \textsc{Lean} \cite{lean} or \textsc{Coq} \cite{coq}, allow users to construct proofs step-by-step whilst the program continuously provides feedback on the logical validity of these steps. For example, if a user tries to show some statement involving a variable which has not been introduced, the system will throw an error. Thus, proof assistants could serve as a training environment for students to freely explore which actions at various stages of a proof are logically allowed. 
Additionally, they actively track the available variables, hypotheses, and proof objectives, which can help students that are unsure about what statements they can use or what needs to be shown.
For some accounts from teachers who have used proof assistants in their courses with the explicit aim to get students to produce better handwritten proofs, see the accounts by Nipkow \cite{Nipkow-VMCAI12}, Patrick Massot \cite[\S 2.2]{french-experience-proof-assistants}, and Fr\'ed\'eric Le~Roux~\cite[\S 2.3]{french-experience-proof-assistants}; Patrick Massot and Heather MacBeth also discuss their personal experiences with using \textsc{Lean} to teach proof writing in the context of mathematical analysis on a panel \emph{Teaching with proof assistants}\footnote{Panel on teaching with proof assistants. Lean Together 2021. Panelists: Jasmin Blanchette, Jeremy Avigad, Julien Narboux, Heather Macbeth, Gihan Marasingha \& Patrick Massot. Available at \url{https://leanprover-community.github.io/lt2021/schedule.html.}} at the Lean Together 2021 meeting.

Despite their potential benefits, there are some issues with the current generation of proof assistants that hinder their integration into conventional proof writing courses.

\begin{enumerate}[(I)]
    \item
    \textbf{Being able to write proofs in a proof assistant does not imply being able to write good-quality proofs by hand.}~~
    Although Thoma and Iannone \cite{Thoma2021} found that students who partook in a workshop on \textsc{Lean} produced better handwritten proofs, a study by Knobelsdorf et al. \cite{Bohne-evidence} found that students in a small course involving \textsc{Coq} performed worse at writing proofs with pen and paper than with \textsc{Coq} itself. The members of the \textsc{Lean} panel attest to a similar statement: students were able to master creating proofs with the proof assistant, but their handwritten proofs left much to be desired.
    Knobelsdorf et al. suggest that, in their case,  the transfer of proof skills might have failed because \textsc{Coq} provides additional scaffolding, like the automated bookkeeping overview, which continuously displays the available assumptions and current proof goals, that was not carefully dismantled during the course.
    
    \item
    \textbf{Proof assistants have a steep learning curve.}~~
    It takes time to learn the specific proof language used by a proof assistant. 
    B\"ohne and Kreitz \cite{Bohne-three-steps}, for example, remark that \textsc{Coq}'s syntax is difficult to learn for students due to the unstructured naming convention of the keywords that indicate certain proof steps.
    Additionally, the foundational system used by most proof assistants, namely type theory, differs enough from the usual set theory to require some explanation. Propositions and subsets in particular are treated differently.

    \item
    \textbf{The feedback provided by proof assistants is reactive and limited to isolated proof steps.}~~ 
    Students need other kinds of feedback as well:
    students might get stuck and need a hint in order to continue with a proof; other students might produce technically correct proofs that are too complicated.
    Think of students neglecting to use lemmas and trying to derive everything from first principles.

    \item
    \textbf{Installing a proof assistant is difficult.}~~The installation instructions can be involved, especially for non-\textsc{Linux} platforms. They often require users to use the command line, but most mathematics freshmen have never used the command line before.
    At our institution, students use personal laptops instead of fixed computers in a computer lab, so pre-installing the proof assistants is not an option.
    Luckily, the developers behind most proof assistants are realizing the importance of having a simple installation procedure, see for example the development of \textsc{Coq Platform}\footnote{Available at \href{https://github.com/coq/platform}{https://github.com/coq/platform}.} which allows users to just install the binaries using a graphical installer.
    
    \item
    \textbf{The user interfaces for proof assistants can be uninviting.}~~Some editors, like the CoqIDE have a dated visual design, giving the impression that the program has not seen maintenance in a while.
    More modern editors exist, both \textsc{Coq} and \textsc{Lean} provide extensions for \textsc{VS Code}, but these interfaces strongly resemble coding environments. Although this might provide computer science students with a sense of familiarity, it can have the opposite effect on mathematics students.

    \item
    \textbf{Many proof assistants do not allow for \LaTeX-formatted expressions.}~~Although proof assistants like \textsc{Coq} and \textsc{Lean} allow for Unicode notation, this is a large downgrade from the beauty and versatility of \LaTeX. 
    This is not a superficial requirement either: good notation clarifies mathematical concepts and aids understanding.
\end{enumerate}

\begin{figure}[htp]
    \centering
    \includegraphics[width=\textwidth]{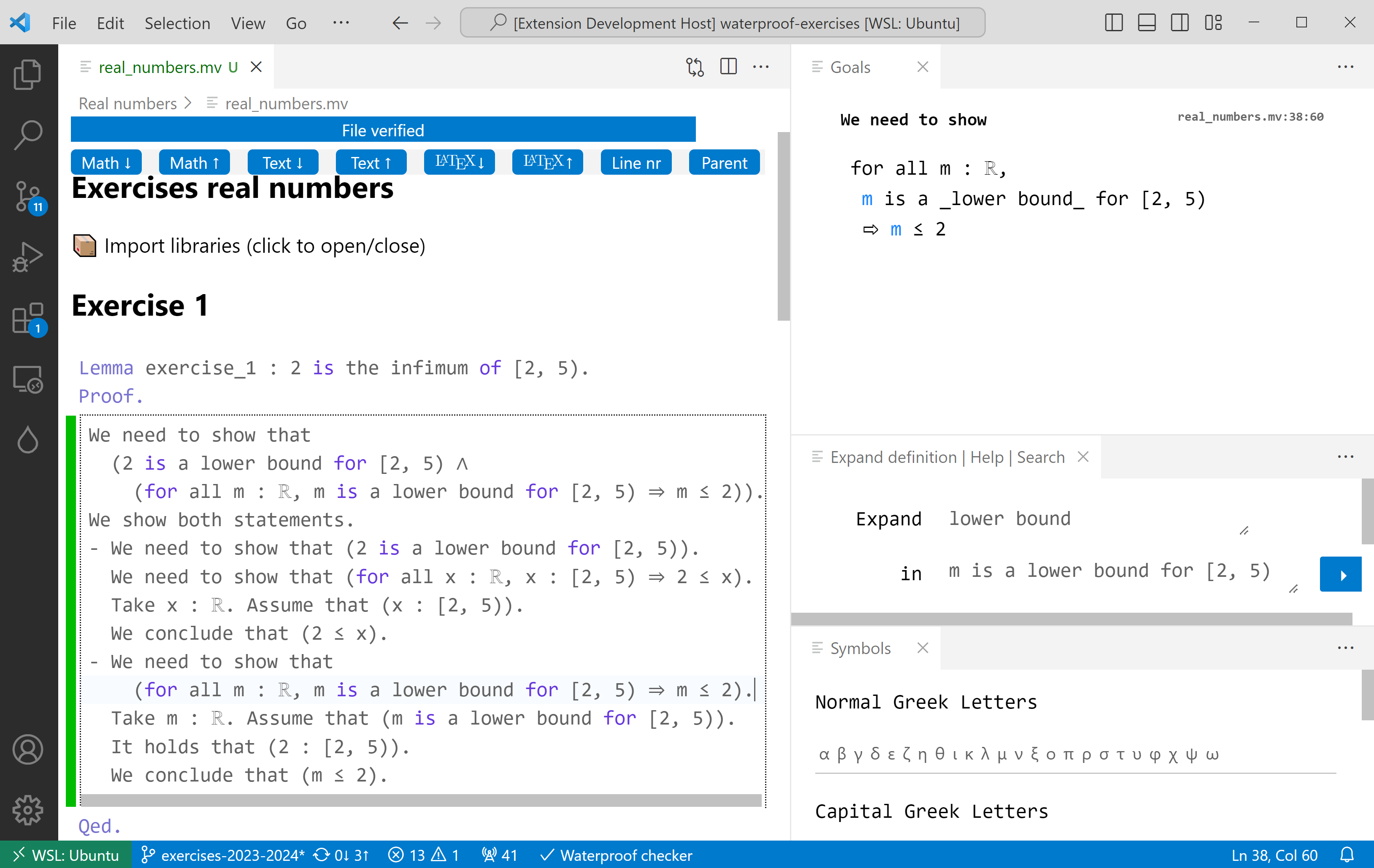}
    \caption{Screenshot of the Waterproof editor showcasing Waterproof's custom proof language. The mixed document with formatted text and verified mathematics is shown on the left. The line prefixed by the box-emoji hides the document preamble that imports the required libraries. The vertical green bar indicates an input area where students can write their proof. The green color indicates that the proof is correct. On the right are multiple panels, from top to bottom they are: the automated bookkeeping overview, limited to only show the proof goal; a panel for expanding definitions, in other editors this is done within the main document; an overview of mathematical symbols, which can be inserted into the main document by clicking on them.}
    \label{fig:screenshot}
\end{figure}

To address the issues above, we created \href{https://impermeable.github.io}{\emph{Waterproof}}, an adaptation of the \textsc{Coq} proof assistant that is designed specifically for teaching mathematical proof writing.
A screenshot is shown in Figure \ref{fig:screenshot}.
We felt that writing a proof in the existing proof assistants differs too much from the ordinary way of writing a proof, that this gap might explain the difficulty with transferring proof skills from proof assistants to pen-and-paper proofs (issue I), and that it and prevents student from tapping into prior knowledge when learning how to use a proof assistant (issue II). Hence, a key idea in Waterproof's design was that
\begin{quote}
    \emph{Writing a proof in Waterproof should be as close as possible to writing a proof by hand, both in terms of the final product and the process of constructing the proof.}
\end{quote}

The Waterproof software\footnote{Available at \href{https://impermeable.github.io}{https://impermeable.github.io}. The custom proof language is part of the \emph{coq-waterproof} plugin, available at \href{https://github.com/impermeable/coq-waterproof}{github.com/impermeable/coq-waterproof}, this article pertains to version 2.1.0+8.17\,. The custom editor is a \textsc{VS Code} extension, available at \href{https://marketplace.visualstudio.com/items?itemName=waterproof-tue.waterproof}{marketplace.visualstudio.com/items?itemName=waterproof-tue.waterproof}, this article pertains to version~1.0.0\,. The editor can be installed directly from \textsc{VS Code}'s `Extensions' panel, after which a walkthrough will guide the user on how to install the coq-waterproof plugin using a graphical installer.} consists of two components: a custom proof language and a custom editor.
The custom proof language allows \textsc{Coq} proofs to be written in a style that closely resembles the style of handwritten proofs. It is part of a custom \textsc{Coq} plugin that also contains subroutines that can give basic suggestions for the next proof step based on the main connector of the proof goal (first steps towards solving issue III).
The custom editor provides a modern user interface (issue V) and uses \emph{mixed documents} which beside verified \textsc{Coq} proofs also contain formatted text (including \LaTeX-expressions (issue VI)). Both the custom editor and proof language are easy to install (issue IV).
Some of the editor's features were designed around its use as an educational tool. Student input, for example, can be limited to specific parts of the document to prevent the accidental deletion of exercises.
Features like an autocomplete function for mathematical symbols and proof steps also help to improve the usability of the editor in general.
The amount of information shown in the automated bookkeeping panel is, by default, also restricted in order to make the proof construction process with Waterproof more closely resemble the process of constructing a proof with pen-and-paper.

Waterproof has been used to supplement the teaching of mathematical analysis at Eindhoven University of Technology (TU/e) for the last four years.
Roughly 175 students register for the course, Analysis~1, every year; the majority are first-year undergraduate students. A selection of the weekly homework exercises have been made available as Waterproof exercise sheets, i.e. documents where the proofs are left to be filled out.
Students can choose to hand in homework by hand or using Waterproof. The homework assignments written in Waterproof are graded automatically.
Use of Waterproof has been optional as only a couple of instructors could answer questions about the program, naturally limiting its adoption by students.
At the start of the 2022-2023 course, 25 student groups ($\approx$ 100 students) handed in homework written in Waterproof; 19 groups ($\approx$ 76 students) continued using it for the final assignment.
We observed that students started using the specific formulations from Waterproof's proof language in their handwritten proofs; the explicit phrasing of these sentences helped to clarify the logical structure of their arguments.

Section~\ref{sec:related-work} explores alternative solutions proposed by others to the problems I--VI outlined above.
Waterproof's custom proof language and its implementation are discussed in Section~\ref{sec:proof-language}, followed by a discussion of the custom editor in Section~\ref{sec:editor}. Section~\ref{sec:use-in-education} outlines how Waterproof was employed in the Analysis~1 course at the TU/e, and reports on  both students' and teachers experience of using the software.

\section{Related Work}
\label{sec:related-work}

As proof assistants have gained popularity among mathematicians, these programs have also found their way into mathematics education. The introduction already mentioned MacBeth and Massot who both used \textsc{Lean} to teach mathematical analysis, as well as the others who spoke at the \textsc{Lean} panel on teaching with proof assistants.
The use of proof assistants in education is not a new phenomenon: communities who actively use these tools in their own research seem to have repeatedly attempted to harness their potential benefits for teaching purposes.
The oldest example we found was the use of Mizar~\cite{mizar} to teach propositional logic at the University of Warsaw in the 70s, according to \cite{mizar-history}.
In computer science, proof assistants are often used to teach more advanced theoretical subjects like formal logic \cite{Martin-Aquinas,Nipkow-VMCAI12} or type theory \cite[\S 3]{Femke-and-Wiedijk}; students report this makes the courses feel more practical, as writing proofs in a proof assistants feels similar to coding.
For some, proof assistants are merely a means to teach the theoretical contents of a course, whereas Nipkow, Massot, Le Roux, and MacBeth also used proof assistants to improve the quality of students' handwritten proofs.

Although some teachers have used existing proof assistants directly, many others saw the need to adapt these programs. Below is a list of tools and techniques that have been developed which would address (parts of) the problems I--VI discussed in the introduction.

B\"ohne and Kreitz (coauthors of the small \textsc{Coq} study \cite{Bohne-evidence} mentioned in Issue I) have developed a didactic method --- opposed to a software solution --- to explicitly guide the transition from \textsc{Coq} proofs to handwritten ones \cite{Bohne-three-steps}.
The gap between formal \textsc{Coq} proofs and informal textbook-style proofs is bridged by introducing three intermediate proof styles that gradually lower the level of formality.
Students are first taught to write proofs in \textsc{Coq}, with custom, easier to learn tactics, and then to both translate between the different levels of formality, as well as develop proofs from scratch at each level.
The authors mention that a solution like Waterproof, which modifies \textsc{Coq} itself to make the proofs be more like handwritten proofs, would be very interesting, but that such a development would require a large amount of preparatory work.
The teaching method, on the other hand, is cheap and flexible: it is easy to adapt to different domains and can be adjusted as a course is being taught.

\textsc{Lurch} \cite{Lurch} markets itself as a word processor that offers proof verification as an additional service, like a spellchecker. It is a total solution with its own deduction engine.
The system does not force a specific syntax on the users, instead semantic information is obtained by having users manually annotate mathematical expressions as `claims', `reasons' or `premises'.
The mathematical expressions can be written using \LaTeX, which is dynamically formatted.
As of now, proof checking is limited to logic and set theory; proofs have to be written out in full detail, as by design \textsc{Lurch} has no automation facilities.
Earlier versions included a computer algebra system for verification as well, but this has been removed due to insufficient precision.

The \textsc{Diproche} system \cite{diproche-1,diproche-2} is similar to \textsc{Lurch}: users can write out proofs in controlled natural language (German); the program parses the text and checks the proof using several of its own deduction engines.
Unique to \textsc{Diproche} is its ability to point out common mistakes and, in some cases, produce counterexamples.
Currently, the \textsc{Diproche} system supports exercises in propositional logic, set theory, elementary number theory, axiomatic geometry and elementary group theory; the parsers and deduction systems had to be adjusted to each domain.
\textsc{Diproche} runs on a remote web server, avoiding the problem of installation.

\textsc{Edukera}\footnote{Available at \href{https://edukera.com}{https://edukera.com}.} \cite{edukera} is a commercial web-application based on the \textsc{Coq} proof assistant that allows users to advance a proof using a click-based graphical interface. Pierre Guillot and Julien Narboux \cite[\S 2.4]{french-experience-proof-assistants} and Simon Modeste \cite[\S 2.5]{french-experience-proof-assistants} describe their experiences of using \textsc{Edukera} in their courses. They note that, due to its point-and-click interface, students are not required to memorize definitions and that some students managed to finish exercises without really understanding what they had to do \cite[\S 2.4]{french-experience-proof-assistants}. At the \textsc{Lean} panel on teaching with proof assistants, Patrick Massot mentioned that he refrained from using \textsc{Edukera} in the past because even he himself started randomly pressing buttons out of frustration with the tool. Guillot and Narboux also mention that \textsc{Edukera} has not been maintained since 2020 and that teachers cannot add their create new exercises on their own.
\vspace{2em}

\textsc{d$\exists\forall$duction}\footnote{Developed by Fr\'ed\'eric Le Roux, Marguerite Bin, Florian Dupeyron \& Antoine Leudi\`ere~(2020). Available at \href{https://perso.imj-prg.fr/frederic-leroux/d\%E2\%88\%83\%E2\%88\%80duction/}{https://perso.imj-prg.fr/frederic-leroux/d\%E2\%88\%83\%E2\%88\%80duction}.} \cite[\S 2.3]{french-experience-proof-assistants} is another click-based graphical interface for constructing proof terms step by step. It is built on top of the \textsc{Lean} theorem prover and runs locally. Exercises can be created using \textsc{Lean} itself, a custom parser allows teachers to specify per exercise which inference rules and definitions are available in the graphical interface. 
\textsc{d$\exists\forall$duction} is not meant as a replacement for conventional exercises, it is a tool that students can use when writing proofs by hand, e.g. to check that the steps they wrote down are allowed or that they have the intended effects.

\textsc{ProofWeb} \cite{Femke-and-Wiedijk} provides web access to a \textsc{Coq} proof assistant running on a central server and adds some additional features tailored to education.
The user interface is similar to that of \textsc{Coq} itself, but the proof state can be displayed in different ways, like Gentzen's deduction trees or Fitch's flag-style proofs.
The central server not only provides an easy way to access \textsc{Coq}, but also serves as a distribution point for exercises and allows teachers to keep track of students' progress.
A parser checks whether students are cheating by using \textsc{Coq}'s powerful automation procedures instead of the intended tactics; these tactics use a custom syntax that is easier to learn, like with B\"ohne and Kreitz in \cite{Bohne-three-steps}.
\textsc{ProofWeb} use \textsc{Coq} version 8.2, which dates to 2009.

\textsc{jsCoq} \cite{jsCoq} is an adaptation of \textsc{Coq} that runs purely in the browser using JavaScript. Since \textsc{jsCoq} requires no installation, it is often used in workshops to introduce people to \textsc{Coq}. \textsc{jsCoq} uses mixed documents that combine executable \textsc{Coq} code with formatted text.
The \textsc{jsCoq} version of the Software Foundations series\footnote{Emilio Jes\'us Gallego Arias, Beno\^it Pin \& Pierre Jouvelot (2017). Available at \href{https://jscoq.github.io/ext/sf}{https://jscoq.github.io/ext/sf}.}, for example, allows students to directly see the effects of executing the \textsc{Coq} code explained in these books.

\textsc{Lean verbose}\footnote{Developed by Patrick Massot (2021). \textsc{Lean 3} version available at \href{https://github.com/PatrickMassot/lean-verbose}{https://github.com/PatrickMassot/lean-verbose}.} \cite[\S 2.2]{french-experience-proof-assistants} is a custom proof language for the \textsc{Lean} theorem prover developed by Patrick Massot that, like Waterproof, mimics natural language formulations for tactics. Its goal is to make it easier for students to transition from proof assistants to pen-and-paper proofs. Like Waterproof (see Section 3.1), \textsc{Lean verbose} also implements a tactic that is able to provide students with suggestions on how to proceed if they are stuck. \textsc{Lean verbose}'s help functionality is more advanced than the one in Waterproof: it is able to offer suggestions on how to use hypotheses and course-specific hints like how to expand certain definitions.

Comparing Waterproof to the solutions above, \textsc{Lurch} and \textsc{Diproche} seem to be most similar: they provide their own editors and allow users to write proofs in a natural language. Both systems, however, use their own proof checkers, hence, they miss out on efforts by the larger community interested in formalizing mathematics. \textsc{Lean verbose} is similar to Waterproof's custom proof language, and was even developed with the same goal in mind, a comparison between the two systems should be interesting. 
B\"ohne and Kreitz provide the most radically different solution: instead of adapting proof assistants to fit education, they suggest a new method to incorporate the existing software into mathematics classes, one that is different from just having students write proofs in both the program and with pen-and-paper. \textsc{d$\exists\forall$duction} is also interesting in this regard since it serves as a tool that students can use to check their reasoning when writing conventional proofs. Like Waterproof, \textsc{jsCoq} allows for the combination of formatted text with \textsc{Coq} proofs in mixed documents; Waterproof was partially inspired by \textsc{jsCoq} and we might try to use \textsc{jsCoq} as a basis for a future version of Waterproof that runs completely in the browser.

\section{Custom Proof Language}
\label{sec:proof-language}

Waterproof's custom proof language allows proofs to be written in a style that closely resembles handwritten proofs. To see how well the handwritten style is approximated, see the example proof shown in Figure \ref{fig:coq-waterproof-proof-text}. Figure \ref{fig:coq-proof-text} shows the same proof written in the default \textsc{Coq} language. By closing the gap between Waterproof's proof style and that of ordinary proofs, we expect that proof skills learned in Waterproof are transferred to pen-and-paper more easily (issue I). The Waterproof language is also easy to learn (issue II): for the Analysis 1 course, the entire language is explained in a single tutorial file, which takes the students only a couple of hours to complete. The custom proof language is part of the coq-waterproof plugin, which also contains the mathematical library used for the Analysis 1 course (a custom library based on the \textsc{Coq} standard library). 

In this section, we often use the term \emph{tactic}. This is proof-assistant terminology for a function that alters the state of the proof, like the subroutine underlying the proof step for introducing a new variable.

\begin{figure}[htp]
    \begin{coqcode}
    Lemma example\_coq\_waterproof :\\
    \phantom{~~}for all $\varepsilon$~:~ℝ, $\varepsilon$ > 0 ⇒ there exists a~:~ℝ, a~:~[0,4) ∧ 4 - $\varepsilon$ < a.\\
    Proof.\\
    \phantom{~~}Take $\varepsilon$~:~ℝ.~Assume that ($\varepsilon$ > 0).\\
    \phantom{~~}Either ($\varepsilon$ < 2)~or ($\varepsilon$ ≥ 2).\\
    \phantom{~~}- Case ($\varepsilon$ < 2).\\
    \phantom{~~-~}Choose a := (4 - $\varepsilon$/2).\\
    \phantom{~~-~}We show both (a~:~[0,4)) and (4 - $\varepsilon$ < a).\\
    \phantom{~~-~}+ We need to show that (0 ≤ a ∧ a < 4).\\
    \phantom{~~-~+~}We show both (0 ≤ a) and (a < 4).\\
    \phantom{~~-~+~}* We conclude that (\& 0 < 4 - 1 < 4 - $\varepsilon$/2 = a).\\
    \phantom{~~-~+~}* We conclude that (a < 4).\\
    \phantom{~~-~}+ We conclude that (4 - $\varepsilon$ < a).\\
    \phantom{~~}- Case ($\varepsilon$ ≥ 2).\\
    \phantom{~~-~}Choose a := 3.\\
    \phantom{~~-~}We show both (3~:~[0,4)) and (4 - $\varepsilon$ < 3).\\
    \phantom{~~-~}+ We conclude that (3~:~[0,4)).\\
    \phantom{~~-~}+ We conclude that (\& 4 - $\varepsilon$ ≤ 4 - 2 = 2 < 3).\\
    Qed.
    \end{coqcode}
    \caption{Proof of $\forall\,\varepsilon > 0~\exists\,a \in [0,4),\, 4 - \varepsilon < a$ written using Waterproof's custom proof language.}
    \label{fig:coq-waterproof-proof-text}

    \bigskip

    \begin{coqcode}
    Lemma example\_coq :\\
    \phantom{~~}forall $\varepsilon$~:~R, $\varepsilon$ > 0 -> exists a~:~R, ([0,4) a) /\char`\\~4 - $\varepsilon$ < a.\\
    Proof.\\
    \phantom{~~}intro $\varepsilon$.~intro $\varepsilon$\_gt\_0.\\
    \phantom{~~}assert ($\varepsilon$ < 2 \char`\\/ 2 <= $\varepsilon$) as cases by lra.\\
    \phantom{~~}destruct cases as [$\varepsilon$\_lt\_two | two\_le\_$\varepsilon$].\\
    \phantom{~~}- \textit{(* Case $\varepsilon$ < 2.~*)}\\
    \phantom{~~-~}exists (4 - $\varepsilon$/2).\\
    \phantom{~~-~}split.\\
    \phantom{~~-~}+ split.\\
    \phantom{~~-~+~}* assert (0 <= 4 - $\varepsilon$/2) as h1 by lra; exact h1.\\
    \phantom{~~-~+~}* assert (4 - $\varepsilon$/2 < 4) as h2 by lra; exact h2.\\
    \phantom{~~-~}+ assert (4 - $\varepsilon$ < 4 - $\varepsilon$/2) as h3 by lra; exact h3.\\
    \phantom{~~}- \textit{(* Case $\varepsilon$ $\geq$  2.~*)}\\
    \phantom{~~-~}exists 3.\\
    \phantom{~~-~}split.\\
    \phantom{~~-~}+ assert (0 <= 3 < 4) as h4 by lra; exact h4.\\
    \phantom{~~-~}+ assert (4 - $\varepsilon$ < 3) as h5 by lra; exact h5.\\
    Qed.
    \end{coqcode}
    \caption{Proof of the same statement written using \textsc{Coq}'s default proof language (written in such a way that shows to what extent a pen-and-paper proof style can be replicated).}
    \label{fig:coq-proof-text}
\end{figure}

\subsection{Features} \label{sec:proof-lang-features}
The main features of Waterproof's custom proof language are listed below, most of these are exemplified by the proof in Figure~\ref{fig:coq-waterproof-proof-text}. They are to be contrasted with the default \textsc{Coq} language in Figure \ref{fig:coq-proof-text}.

\begin{itemize}
    \item \textbf{Proof step formulations inspired by handwritten proofs.}~~
    The proof steps formulations in Waterproof's custom proof language take the form of full sentences that are used in ordinary mathematical proofs as well.
    Compare the formulations used by Waterproof's language~(left) and \textsc{Coq}'s default \textsc{Coq} syntax (right) for introducing the $\varepsilon$ variable:
    \vspace{2pt}
    \begin{center}
    \begin{tblr}{width=0.9\linewidth,colspec={X[3,l]X[1,c]X[3,l]}}
        \SetCell{gray!20}\texttt{Take $\varepsilon$~:~$\mathbb{R}$.} & vs. & \SetCell{gray!20}\texttt{intro $\varepsilon$.}
    \end{tblr}
    \end{center}
    \vspace{3pt}
    The Waterproof formulation includes the mathematically relevant information that $\varepsilon$ is a real number, whereas the default \textsc{Coq} tactic does not. The difference becomes even more apparent if we compare the formulations for introducing the assumption that $\varepsilon > 0$\,:
    \vspace{2pt}
    \begin{center}
    \begin{tblr}{width=0.9\linewidth,colspec={X[3,l]X[1,c]X[3,l]}}
        \SetCell{gray!20}\texttt{Assume that ($\varepsilon$ > 0).} & vs. & \SetCell{gray!20}\texttt{intro $\varepsilon$_gt_0.}
    \end{tblr}
    \end{center}
    \vspace{3pt}
    The actual proposition itself is not even mentioned in the default \textsc{Coq} formulation. For \textsc{Coq} users, this is a non-issue, since the content of the assumption can easily be found in the automated bookkeeping overview. For educational purposes, however, we wish to limit the information shown there such that students learn to write proofs without it. To make the \textsc{Coq} proofs readable on their own, the relevant mathematical information has to be included back into the proof steps themselves. In fact, with Waterproof's proof step formulations, we found ourselves to pay less and less attention to the automated bookkeeping overview anyway.

    \phantom{~~~~~~}
    Like in ordinary proofs, Waterproof uses different formulations for introducing a variable and introducing an assumption. Default \textsc{Coq} uses the same tactic for both, because to its foundational system, type theory, there is no meaningful difference between the two. In Waterproof, an error is thrown if the wrong formulation is used, like for example in
    \begin{coqcode}
        Take $\varepsilon$_gt_0~:~($\varepsilon$ > 0).
    \end{coqcode}
    
    \item \textbf{Implicit use of automation to verify statements.}~~
    Basic statements do not need to be justified in ordinary proofs, but, by design, proof assistants require a proof for every claim that is made. To get Waterproof's proof writing style more in line with that of handwritten proofs, and to prevent students from being bogged down by having to show `obvious' statements, Waterproof's tactics try to prove the user's claims automatically, like in the proof steps
    \begin{coqcode}
        It holds that (...).\\
        We conclude that (...).\\
        It suffices to show that (...).
    \end{coqcode}
    \vspace{-7pt}
    Explicit justifications can still be provided by referring to specific lemmas, for example by writing 
    \begin{coqcode}
        By lemma_1 it holds that (...).
    \end{coqcode}
    \vspace{-6pt}

    \phantom{~~~~~~}
    Using an automation system brings multiple advantages. First, users are able to write proofs without them having to know all the details of the mathematical library. They do not need to know the specific names for basic properties like $0 < 1$\,. We also found that with automation it became much easier to use a forwards-reasoning style, which better matches the flow of regular math proofs than the backwards-reasoning style encouraged by most proof assistants.
    Finally, the automation system allows Waterproof's proof step formulations to de-emphasize the use of labels. In the default \textsc{Coq} language, every statement is labelled, because these labels are needed to construct the explicit justifications required by the proof assistant.
    
    \phantom{~~~~~~}
    In the creation of mathematical libraries and exercises for their courses, teachers can tune the automation system such that it is able to show those (and only those) statements which they consider to be `obvious'. The lemmas and solvers that the system is allowed to use in its proof search are managed using \textsc{Coq}'s so-called \emph{hint databases}, which can be customized by teachers themselves. For the Analysis 1 course, for example, we added custom tactics for verifying computations involving absolute values.

    \phantom{~~~~~~}
    Some statements are actively shielded from the full power of the automation system. Only a selection of the available hint databases is used to try to prove such statements automatically. We chose to shield statements starting with a logical operator, e.g. $\forall\, x : \mathbb{R},\, x^2 \geq 0$\,, since this forces students to actively engage with such operators for non-trivial statements.
    
    \item \textbf{Mandatory signposting.}~~
    Handwritten proofs contain many statements that indicate what needs to be shown or what case is being considered.
    Proofs written in proof assistants often neglect such \emph{signposting} altogether, because the current proof state and goal are constantly displayed in the automated bookkeeping overview.
    Proofs can be signposted using comments, but these are optional. In contrast, Waterproof's proof language at certain points forces users to signpost their proofs, thus we communicate to students that we expect them to use these statements in their handwritten proofs as well.

    \phantom{~~~~~~}
    One example is the mandatory inclusion of a statement indicating what case is being considered after the case distinction $(\varepsilon < 2) \vee (\varepsilon \geq 2)$ in Figures \ref{fig:coq-waterproof-proof-text} and \ref{fig:coq-proof-text}.
    In Waterproof's custom proof language, the \colorbox{gray!20}{\texttt{Case ($\varepsilon$ < 2).}} statement \emph{has} to be included before further progress can be made. The corresponding sentence in the default \textsc{Coq} proof,  \colorbox{gray!20}{\texttt{(*~Case $\varepsilon$ < 2.~*)}}, is merely a comment. Use of the \colorbox{gray!20}{\texttt{Case (...).}} statement in Waterproof is enforced by the \colorbox{gray!20}{\texttt{Either (...)~or (...).}}-tactic which performs the case distinction. It changes the first subgoal, which is visible to the user via the automated bookkeeping overview, to the text
    \begin{coqcode}
        Add the following line to the proof:\\
        \phantom{~~}Case ($\varepsilon$ < 2).
    \end{coqcode}
    \vspace{-7pt}
    After the correct \colorbox{gray!20}{\texttt{Case (...).}}-statement is included, the text is changed back to show the original goal.

    \phantom{~~~~~~} Besides forcing the users to indicate what case is being considered after a case distinction, Waterproof's proof language enforces users to indicate the base case and the induction step after the announcement of an induction argument, and to indicate what goal is being show after the proof of a $\wedge$-statement is split into separate subproofs.
    
    \phantom{~~~~~~} Users can also insert optional signposting by using the statement
    \begin{coqcode}
        We need to show that (...).
    \end{coqcode}
    \vspace{-7pt}
    whenever they feel that this helps to keep the proof readable. This sentence checks whether the goal specified by the user actually corresponds to what needs to be shown, thus acting as sanity-check.

    \item \textbf{More elaborate error messages.}~~
    Waterproof's custom tactics provide elaborate feedback to help students fix issues themselves. The errors thrown by \textsc{Coq}'s default tactics can sometimes be difficult to parse. For example, if the goal is to show that $\exists\, x : \mathbb{R},\, x \geq 0$\,,
    and the default \textsc{Coq} tactic \colorbox{gray!20}{\texttt{intro x.}} is mistakenly used, a cryptic error message is returned:
    \begin{coqcode}
        No product even after head-reduction.
    \end{coqcode}
    \vspace{-7pt}
    In comparison, the equivalent Waterproof tactic \colorbox{gray!20}{\texttt{Take $x$~:~$\mathbb{R}$.}} returns the error message
    \begin{coqcode}
        `Take ...'~can only be used to prove a `for all'-statement ($\forall$) \\
        ~~~~or to construct a map ($\to$).
    \end{coqcode}
    \vspace{-7pt}
    Note that Waterproof's proof steps formulations also allow for more mistakes to be made since they require more data to be specified by the user.     
    For example, if the user starts a $\forall\, x : \mathbb R, ...$ proof with \colorbox{gray!20}{\texttt{Take $x$~:~$\mathbb N$}.}, they are presented with the error
    \begin{coqcode}
        Expected a variable of type ℝ instead of ℕ.
    \end{coqcode}
    \vspace{-7pt}
    Such a mistake would not have been possible to make with \textsc{Coq}'s default `\texttt{intro}'-tactic.
        
    \item \textbf{Conventional mathematical notation.}~~
    Waterproof's custom proof language provides basic notations for common mathematical objects like sets ($\mathbb N$\,, $\mathbb Z$\,, $\mathbb R$, etc.), logical operators ($\forall$\,, $\exists$, etc.), and intervals. For function application, the mathematical convention is used instead of the functional programming-notation used by default \textsc{Coq}, e.g.
    \vspace{2pt}
    \begin{center}
    \begin{tblr}{width=0.9\linewidth,colspec={X[3,l]X[1,c]X[3,l]}}
        \SetCell{gray!20}\texttt{f(1+1,2)} & vs. & \SetCell{gray!20}\texttt{f (1+1) 2}
    \end{tblr}
    \end{center}
    \vspace{3pt}
    We also provide an alternative way to denote subset membership than the one commonly used in default \textsc{Coq} (and other proof assistants based on type theory), namely
    \vspace{2pt}
    \begin{center}
    \begin{tblr}{width=0.9\linewidth,colspec={X[3,l]X[1,c]X[3,l]}}
        \SetCell{gray!20}\texttt{x~:~[0,1)} & vs. & \SetCell{gray!20}\texttt{[0,1) x}
    \end{tblr}
    \end{center}
    \vspace{3pt}
    Both expressions are definitionally equal to the proposition $0 \leq x < 1$\,. We wanted the :-symbol to just read to students as a replacement for $\in$-symbol, without having to explain the logical differences.
    Since some may object to the abuse of the type-theoretical :-symbol to denote subset membership, this notation has been made optional.

    \phantom{~~~~~~}
    The custom language also provides a framework for expanding definitions whose notation consists of multiple words, like `\texttt{1 is the supremum of [0,1)}'. 
    To expand this definition in default \textsc{Coq}, the user has to know that this term is represented internally as `\texttt{is_sup [0,1) 1}'\,, since it is the term `\texttt{is_sup}' that needs to be expanded. With Waterproof's custom proof language it is possible to write
    \begin{coqcode}
        Expand the definition of supremum in (1 is the supremum of [0,1)).
    \end{coqcode}
    \vspace{-7pt}
    which returns the expanded definition. To allow users to unfold notation in this way does require some extra lines of code to be added to mathematical libraries in places where notations are defined, but this is within reach of teachers who are developing their own libraries.

    \item \textbf{Chains of (in)equalities.}~~Chaining simple (in)equalities to prove more complicated ones is key ingredient of mathematical proofs, especially in analysis. Some proof assistants support this style of reasoning, like \textsc{Lean}'s calculational proofs, but \textsc{Coq}, by default, does not. Waterproof's custom proof language does provide the ability to use such chains, like in the sentence
    \begin{coqcode}
        We conclude that (\& 0 < 4 - 1 < 4 - $\varepsilon$/2 = a).
    \end{coqcode}
    \vspace{-7pt}    

    \item \textbf{Suggestion for next proof step.}~~
    In a first attempt to provide more than just reactive feedback, the custom proof language is equipped with a \colorbox{gray!20}{\texttt{Help.}}-tactic that gives basic suggestions for what the next step in a proof should be. The suggestion is based on the current shape of the goal. For example, if the goal is to show that $\forall\, x : \mathbb R,\,x^2  \geq 0$\,, the suggestion will be to introduce the variable~$x$ using the tactic \colorbox{gray!20}{\texttt{Take x~:~$\mathbb R$.}}.

\end{itemize}

\subsection{Implementation}
This section highlights some of the implementation details of Waterproof's custom proof language, parts of the exposition can be quite technical. The custom language is mostly implemented using standard \textsc{Coq} scripts, the \textsc{Ltac2} tactic language is used to specify the behavior of the custom tactics. The automation system is partially written in \textsc{OCaml}, as this allowed for better access to \textsc{Coq}'s internal structure like how the hint databases are called by the default `\texttt{auto}'-tactic.

Our coq-waterproof library provides many examples of custom \textsc{Ltac2} tactics that are non-standard in the sense that they are not-necessarily aimed at manipulating the proof term. For example, large parts of the implementation of the \colorbox{gray!20}{\texttt{Take ...}}-tactic discussed and shown below do not alter the proof state but instead serve to give detailed error messages.
We hope that our library can provide inspiration to people learning the \textsc{Ltac2} tactic language.

\subsubsection{Example of a tactic implementation and the numerous amount of checks on user input}

The \colorbox{gray!20}{\texttt{Take ...}}-tactic is one of Waterproof's natural language tactics implemented using \textsc{Ltac2}, it implementation is shown in Listing \ref{lst:take}. The variables are introduced by calling the default `\texttt{intro}'-tactic (line~7), but before doing so, a numerous amount of checks are performed in order to provide detailed error messages. Thus, the implementation takes up an entire page, instead of a few lines.

The notation for the tactic is defined at the bottom (line 46). It allows multiple variables from different types to be introduced in a single step, e.g. by writing
\begin{coqcode}
    Take n~:~ℕ and x, y~:~ℝ and z~:~ℂ.
\end{coqcode}
\vspace{-7pt}

Per variable and per type, it is checked whether the goal (still) requires a variable of this type to be introduced (lines 2--6, and lines 15--16 together with 22--23). 
For each type, it is also checked whether the type that needs to be introduced is not a proposition (lines 18--19 and 34--37). A type is considered to be a proposition if it belongs to the sort `\texttt{Prop}'. In general, mathematical libraries use `\texttt{Prop}' to encode propositions, but such a neat separation is not always guaranteed: the Homotopy Type Theory library~\cite{hott-lib}, for example, places propositions in the same universe as one-element sets. Note that if the first check (lines 34--37) fails, a different error message is returned than if the repeated check (lines 18--19) fails.

\begin{listing}[htp]
\begin{coqlinenumbers}
   Local Ltac2 intro_ident (id : ident) (type : constr) :=
     lazy_match! goal with
       | [ |- forall _ : ?u, _] =>
         let ct := get_coerced_type type in
         (* Check whether we need a variable of type [type], including coercions. *)
         match check_constr_equal u ct with
           | true  => intro $id
           | false => throw (too_many_of_type_message type)
         end
       | [ |- _] => throw (too_many_of_type_message type)
     end.

   Local Ltac2 intro_per_type (pair : ident list * constr) :=
     let (ids, type) := pair in 
     lazy_match! goal with
       | [ |- forall _ : ?u, _] => 
         (* Check whether [u] is not a proposition. *)
         let sort_u := get_value_of_hyp u in
         match check_constr_equal sort_u constr:(Prop) with
           | false =>
             (* Check whether we need variables of type [type], including coercions. *)
             let ct := get_coerced_type type in
             match check_constr_equal u ct with
               | true  => List.iter (fun id => intro_ident id type) ids
               | false => throw (expected_of_type_instead_of_message u type)
             end
           | true  => throw (of_string "Tried to introduce too many variables.")
         end
       | [ |- _ ] => throw (of_string "Tried to introduce too many variables.")
     end.

   Local Ltac2 take (x : (ident list * constr) list) := 
     lazy_match! goal with
       | [ |- forall _ : ?u, _] => 
         (* Check whether [u] is not a proposition. *)
         let sort_u := get_value_of_hyp u in
         match check_constr_equal sort_u constr:(Prop) with
           | false => List.iter intro_per_type x
           | true  => throw (of_string 
           "`Take ...` cannot be used to prove an implication. Use `Assume that ...` instead.")
         end
       | [ |- _ ] => throw (of_string
       "`Take ...` can only be used to prove a `for all`-statement or to construct a map.")
     end.

   Ltac2 Notation "Take" vars(list1(seq(list1(ident, ","), ":", constr), "and")) := 
     panic_if_goal_wrapped ();
     take vars.
\end{coqlinenumbers}
\caption{Implementation of the \colorbox{gray!20}{\texttt{Take ...}}-tactic used for introducing variables.}
\label{lst:take}
\end{listing}

\subsubsection{Automation system written in \textsc{OCaml}}
The implicit automation system uses a custom version of \textsc{Coq}'s `\texttt{auto}'-tactic that tries to include a user-specified lemma in its proof finding algorithm. Like the original `\texttt{auto}'-tactic, the custom version is written in \textsc{OCaml}. The condition that the user-specified lemma has to be used, allows Waterproof to reject statements like
\begin{coqcode}
    By IVT it holds that (1 + 1 = 2).
\end{coqcode}
\vspace{-6pt}
which claims that $1 + 1 = 2$ holds because of the intermediate value theorem (IVT). The custom `\texttt{auto}'-tactic will try to apply the `\texttt{IVT}'-lemma at every branch of the search tree, and fails if no successful path is found that involves the lemma.

Using \textsc{OCaml} scripts, we are also able better control what hint databases are used by the automation system. In the default \textsc{Coq} language, the databases that are used are specified in the call to the `\texttt{auto}'-tactic directly, but for the Waterproof language to mimic the style of handwritten proofs, this is not possible. Our solution was to create a custom database management system. Based on the context, three different collections of hint databases can be used by the automation system. The `weak' collection, for example, is used for proving shielded statements. We include custom vernacular commands that allow teachers to customize which hint databases are included in these collections when creating their exercises.

The \textsc{OCaml} scripts for the custom `\texttt{auto}'-tactic and the hint database management system were written by Balthazar Patiachvili during his 2023 internship, the technical details can be found in his M1~internship report \cite{balthazar}.

\subsubsection{Enforcing the use of signposting by wrapping the proof goal}

The \colorbox{gray!20}{\texttt{Either (...)~or (...).}}-tactic forces users to indicate what case is being considered by secretly wrapping the proof goal in a type family that only the correct \colorbox{gray!20}{\texttt{Case (...).}}-tactic can undo. The wrapped type is given a print-only notation that informs the user which tactic to use to unwrap the goal. This technique is used by all tactics that enforce the use of signposting. Listing \ref{lst:wrapper} shows the definition of the `case'-wrapper, its print-only definition, and how it is used by (subroutines of) the \colorbox{gray!20}{\texttt{Either (...)~or (...).}} and \colorbox{gray!20}{\texttt{Case (...).}}-tactics to wrap and unwrap the goal.
To make sure that the other tactics in Waterproof's custom proof language do not alter wrapped goals, they first call a `\texttt{panic_if_goal_wrapped}' function, which throws an error if the goal is wrapped. The \colorbox{gray!20}{\texttt{Take ...}}-tactic can be seen to call this function in line 47 of Listing \ref{lst:take}.

\begin{listing}[htp]
\begin{coqlinenumbers}
   Module Case.
     Private Inductive Wrapper (A G : Type) : Type :=
       | wrap : G -> Wrapper A G.
     Definition unwrap (A G : Type) : Wrapper A G -> G :=
       fun x => match x with wrap _ _ y => y end.
   End Case.
 
   Notation "'Add' 'the' 'following' 'line' 'to' 'the' 'proof:' 'Case' ( A )." :=  
     (Case.Wrapper A _) (at level 99, only printing, 
       format "'[ ' Add  the  following  line  to  the  proof: ']' '//'   Case  ( A ).").
\end{coqlinenumbers}
\begin{coqlinenumbersstart}{11}
   Ltac2 either_or_prop (t1:constr) (t2:constr) :=
     let h_id := Fresh.in_goal @_temp in
     let attempt () := assert ($t1 \/ $t2) as $h_id by (* call to automation system *) in
     match Control.case attempt with
       | Val _ =>
         let h_val := Control.hyp h_id in
         destruct $h_val;
         Control.focus 1 1 (fun () => apply (Case.unwrap $t1));
         Control.focus 2 2 (fun () => apply (Case.unwrap $t2))
       | Err exn => throw (of_string
         "Could not find a proof that the first or the second statement holds.")
     end.
\end{coqlinenumbersstart}
\begin{coqlinenumbersstart}{23}
   Ltac2 case (t:constr) :=
     lazy_match! goal with
       | [|- Case.Wrapper ?v _] =>
         match check_constr_equal v t with
           | true => apply (Case.wrap $v)
           | false => throw (of_string "Wrong case specified.")
         end
       | [|- _] => throw (of_string "No need to specify case.")
     end.
\end{coqlinenumbersstart}
    \caption{Implementation of the `case'-wrapper that forces users to state what case is being considered. Top: definition of the wrapper and its print-only notation. Middle: subroutine used by the \colorbox{gray!20}{\texttt{Either (...)~or (...).}}-tactic to wrap the goal by applying `\texttt{Case.unwrap}' in lines 18 and 19. Bottom: subroutine called by the\colorbox{gray!20}{\texttt{Case (...).}}-tactic to unwrap the goal by applying `\texttt{Case.wrap}' in line~27.}
    \label{lst:wrapper}
\end{listing}

\subsubsection{Flexible chains of (in)equalities using typeclasses}

Waterproof's custom notation for chains of (in)equalities makes extensive use of typeclasses.
Chains consists of a base-component, like `\texttt{x < y}', to which multiple link-components, like 
`\texttt{$\leq$ z}' or `\texttt{= w}' can be added. The ordering-symbols in these components are purely formal, they are given a context specific interpretation at a later stage. The `\texttt{\&}'-symbol at the start of an (in)equality chain is short-hand notation for a complicated expression that specifies how base- and link-components are combined using special `\texttt{chain_link}'-operators. 
These components, however, cannot be combined freely: the variables need to be of the same type, and the final combination is not allowed to contain both `\texttt{<}'- and `\texttt{>}'-symbols. This is where typeclasses come in: the `\texttt{chain_link}'-operators are typeclass properties which are only defined for certain combinations of chains and link-components. Typeclass resolution will fail for non-valid combinations. The `\texttt{\&}'-symbol notation also adds a function called `\texttt{total_statement}', which turns a chain, with purely formal ordering symbols, into a large conjunction of individual, semantically meaningful (in)equalities. 
This function is again a typeclass property, such that the interpretation it provides to the formal ordering symbols can easily be adjusted to new contexts by the creators of mathematical libraries.

To function properly, the implementation of (in)equality chains requires some coupling with the rest of the proof language.
In the \colorbox{gray!20}{\texttt{We conclude that (...).}}-tactic, a piece of code is included to check whether the global statement of a chain, for example $0 < a$ for the chain $0 < 4 - 1 < 4 -\varepsilon/2 = a$\,, actually matches the current proof goal.
Secondly, the automation system contains code that first splits chains into their individual components before proving them, so that, if the automation system fails, the user knows exactly which (in)equality is to blame.

\section{Custom Editor}
\label{sec:editor}

The Waterproof editor can be used for \textsc{Coq} documents in general, but it was also designed for use in education specifically.
It uses mixed documents that combine formal \textsc{Coq} proofs with formatted text, including \LaTeX-expressions. The editor offers general quality-of-life improvements, like an autocomplete function for mathematical symbols and common proof steps. For use of documents as exercise sheets, student input can be limited to designated areas. The amount of information shown in the automated bookkeeping panel is by default limited to only showing the proof goal. By removing this scaffolding, we expect that it is easier to transfer proof skills from Waterproof to pen-and-paper proofs (issue I).

\subsection{Features}
The main features of the custom editor are listed below.

\begin{itemize}    
    \item \textbf{Mixed documents.}~~
    The editor was built to support documents that combine verified \textsc{Coq} proofs with easily readable, formatted text. The styling of the text segments is specified using \textsc{Markdown}, which besides basic text formatting also supports bullet points, URLs, images, and animations. \LaTeX-expressions are supported as well, both as inline and displayed equations.

    \item \textbf{Designated input areas.}~~
    In order to use Waterproof documents as \emph{exercise sheets}, documents where the proofs have been left empty to be completed by students, it is important that students' input is limited to designated \emph{input areas}. Otherwise, students might accidentally change or even delete part of the exercises.
    An example of such an input areas is shown in Figure \ref{fig:screenshot}, it is demarcated by the green bar to its left. The green color indicates that the proof provided there is correct; if the proof were to contain an error, the bar would be red. 

    \phantom{~~~~~~} Editing the document outside of designated input areas can be toggled with a setting called `Teacher Mode', it is disabled by default. In principle, students could enable Teacher Mode in the settings menu, but they do not know about its existence, nor would it give them any benefit in our course setup. Instead of allowing Teacher Mode to be set by all users, in a future version we would like to give control over this setting to teachers only.
    
    \phantom{~~~~~~}For the creation of  exercise sheets by teachers, it is convenient to first create a master file that tests the feasibility of exercises. This file can be explicitly saved as an exercise sheet, which removes all the content from input areas.
    
    \item \textbf{Hidden segments.}~~
    Parts of a document can be marked as \emph{hidden segments}, their contents can be revealed and re-hidden by clicking on them. Originally this feature was meant for including optional hints for exercises, but it turned out to also be useful for hiding the preambles that import libraries at the start of a file (for an example, see Figure \ref{fig:screenshot}). If these preambles are not hidden, they distract and confuse students.
    
    \item \textbf{Limited automated bookkeeping.}~~
    Proof assistants alleviate the mental workload for mathematicians by keeping track of variables, assumptions, intermediate assertions, and proof goals in an automated bookkeeping overview. Students, however, might learn to become dependent on these features \cite{Bohne-evidence}. By default, the Waterproof editor only shows the proof goals, the remaining information can be found in a separate panel which needs to be opened on purpose. Whereas with the default \textsc{Coq} language, the automated bookkeeping panel may contain crucial information, Waterproof's custom language itself requires all the information needed to complete a proof.
    
    \item \textbf{Quick and easy input methods for mathematical symbols and proof steps.}~~
    The Waterproof editor provides an autocompletion functionality for mathematical symbols and proof steps from the custom proof language. To type the $\infty$-Unicode symbol, just start typing the string \texttt{\textbackslash infty} and select the correct suggestion from the autocomplete menu.
    Proof steps are even easier thanks to the use of \emph{code snippets}. Code snippets are templates with gaps for context-specific information, the \texttt{tab}-key can be used to move directly between these holes. Autocompletion for code snippets works the same way as for normal strings.
    
    \phantom{~~~~~~}
    Symbols and proof steps can also be inserted by selecting them from one of the additional side panels provided by the editor, see Figure \ref{fig:screenshot}. Although using the autocomplete function is faster, beginning users will not be familiar with the available proof steps and the (\LaTeX) encodings used for the mathematical symbols. These panels provide an overview of which symbols and proof steps are available, they also include brief explanations and examples of how the proof steps are used.

    \item \textbf{Suggestion for next proof step.}~~
    Right-clicking in the proof makes a small `Help' button appear, clicking this button calls the corresponding-tactic from Waterproof's custom proof language (Section \ref{sec:proof-lang-features}, final bullet point), which provides the user with a suggestion what proof step to perform next.

    \item \textbf{Separate panel for expanding definitions.}~~
    The editor provides a way to expand definitions that is more similar to how this is done with pen-and-paper. In most proof assistants, expanding a definition is treated akin to other proof steps: a command in the proof text instructs the program to expand the definition and the result is shown in the information panel tracking the current proof state. When writing a proof by hand, looking up definitions in e.g. the lecture notes, is a separate activity that is spatially removed from the place where the proof is being written. The Waterproof editor provides a dedicated side-panel (see Figure \ref{fig:screenshot}), away from the main proof writing window, where users can instruct the proof assistant to expand definitions and inspect the results.
    
    \item \textbf{Modern, friendly design.}~~
    The Waterproof editor adheres to \textsc{VS Code}'s minimalist UI design and is compatible with its different color themes. Although \textsc{VS Code} is often used for coding, thanks to the use of mixed documents the Waterproof editor does not resemble a coding environment.

    \item \textbf{Compatibility with the \textsc{Coq} ecosystem.}~~
    The Waterproof editor uses \textsc{Coq}'s new \texttt{.mv}-files which support the mixing of formatted text and formal proofs. The editor does not yet support the older, more common \texttt{.v}-files.

    \item \textbf{Optional use of the custom proof language.}~~
    Although by default the editor is configured to use Waterproof's custom proof language, this can be changed in the settings menu. Doing so changes the autocomplete functionality and the proof step-overview panel to show a selection of \textsc{Coq}'s default formulation of proof steps. Without the plugin implementing the custom proof language, the editor is no longer able to provide suggestions for the next proof step, but the separate panel for expanding definitions still works.

\end{itemize}

\subsection{Implementation}

The Waterproof editor uses \textsc{Coq}'s new \texttt{.mv}-files as underlying documents. These files natively combine formatted text and formal \textsc{Coq} proofs. The formatting of the text is specified using the  \textsc{Markdown} language\footnote{Developed by John Gruber (2004). Available at \href{https://daringfireball.net/projects/markdown}{https://daringfireball.net/projects/markdown}.}, HTML-like markers are added to the text parts to indicate the input areas and hidden segments. For example, an input area is marked by \texttt{<input-area>~...~<\textbackslash input-area>}\,. Editing and rendering of the text and code parts respectively is done using the \emph{What You See Is What You Get (WYSIWYG)}-style editors \textsc{ProseMirror}\footnote{Developed by Marijn Haverbeke. Available at \href{https://prosemirror.net}{https://prosemirror.net}.} and \textsc{CodeMirror}\footnote{Developed by Marijn Haverbeke. Available at \href{https://codemirror.net}{https://codemirror.net}.}. The \textsc{prosemirror-math} package\footnote{Developed by Benjamin Bray. Available at \href{https://benrbray.com/prosemirror-math/}{https://benrbray.com/prosemirror-math}.} is used to render \LaTeX-expressions.

The editor uses \textsc{Coq LSP}\footnote{Developed by Emilio Jes\'us Gallego Arias (2023). Available at \href{https://github.com/ejgallego/coq-lsp}{https://github.com/ejgallego/coq-lsp}.}, developed by Emilio Jes\'us Gallego Arias, to communicate with a \textsc{Coq} instance running in the background. \textsc{Coq LSP} implements an extended version of the standard Language Server Protocol (LSP), the extension for example also specifies how the information in the automated bookkeeping panel is communicated between the editor and language server. \textsc{Coq LSP} succeeds \textsc{Coq SerAPI} \cite{serapi}, a \textsc{Coq}-specific communication protocol, which was used by previous versions of the Waterproof editor.
Among other things, \textsc{Coq LSP} supports the new \texttt{.mv}-files and allows for continuous checking of documents. Before, the program had to be instructed manually to verify proof steps, which had to be explained to students. \textsc{Coq LSP} also provides custom ways to deal with errors. For example, incomplete proofs are automatically interpreted as being `admitted', meaning that users can continue working on other proofs without having to explicitly note that they (temporarily) gave up on a preceding proof.

The editor is implemented as a custom \textsc{VS Code} extension since this makes it easier for our small development team to maintain.
Previous versions of the editor were stand-alone programs, which meant that we had to write and maintain code for all kinds of basic functionalities, like menu bars. Such features are now inherited from the \textsc{VS Code} editor, and we can also reuse the standard \textsc{VS Code} language client for communication with \textsc{Coq LSP}. Delegating the maintenance of these code bases allows us to focus on features specific to Waterproof.

We tried to keep the design of the editor as modular as possible to facilitate reuse of its code base in different experiments with editors for proof assistants. Similarly, we have also tried to reduce the coupling between the Waterproof editor and the custom proof language, but some dependency remains: the configuration of the autocompletion and the overview side-panel have to be manually adjusted to match the formulations in Waterproof's proof language. 
The feature that suggest the next proof step also relies on the custom proof language for its implementation.

\section{Use in Education}
\label{sec:use-in-education}
For the last four years at Eindhoven University of Technology (TU/e), Waterproof has been used to supplement teaching the Analysis 1 course.
The goal is to teach students how to rigorously prove theorems from calculus. This includes teaching them about mathematical concepts like metric spaces, sequences, series, convergence, limits, and continuity, but also teaching students how to write proper mathematical proofs. Students have encountered proofs and logic in an earlier course, but their proof writing skills are still in development. This section outlines how Waterproof is used in Analysis 1, how students and teachers experience working with the software, and a brief observation of the effects on students' handwritten proofs.

\subsection{Course Details}

Analysis 1 is a mandatory course for all first-year mathematics undergraduate students at the TU/e; each year approximately 175 students register for the course. Most are first-year students, some are taking the course for a second time. Almost all students are mathematics majors, some do a combined program with either physics or computer science. The course lasts for eight weeks (excluding the exam weeks) and has a study load of 5 ETCS. It consists of biweekly lectures and instruction classes, an exam, a midterm, and weekly homework exercises. The homework exercises are handed-in in groups of four and make up 10\% of the final grade.

\subsection{Waterproof in Analysis 1}

Waterproof is used as an alternative for the regular homework exercises. A selection of exercises have been made available as Waterproof documents where the proofs are left incomplete. Instead of handing in handwritten solutions, students can choose to hand in completed versions of these documents.

Use of Waterproof is not mandatory, and only some of the instructors are able to offer support for Waterproof. 
Last year (2022-2023), three out of the six tutors were able to help students this way.
Not all instructors were intimately familiar with the software: they had only used Waterproof to solve the exercises themselves, but this did allow them to answer many of the students' questions. In case these instructors could not figure out how to solve an issue, they could defer to some who had more experience with Waterproof.
At the start of the course, students are informed which instruction classes will support the use of Waterproof, and that they should register for these classes if they wish to use the program; in the other instruction classes, students are expected to hand in their homework the regular way, i.e. handwritten or typed up in \LaTeX.

No additional time is reserved to teach students how to use Waterproof: they are provided with a tutorial file (which takes a couple of hours to complete) that teaches them how to use Waterproof's custom proof language, and videos explaining how to use the editor. They can also ask questions during their scheduled instruction classes.

Homework solutions handed in as Waterproof files are graded automatically using \textsc{Momotor}, an automated processing tool developed by the TU/e, that is integrated into the \textsc{Canvas} Learning Management System (LMS). Grading is done on a binary scale: if a proof is deemed correct by an instance of the proof assistant running on the \textsc{Momotor} server, full points are awarded; zero points are awarded if the proof contains a mistake or is incomplete. This strict grading scheme matches the feedback provided by the program to students when working on a proof; with human graders there is more room for nuance.
To make sure that a student did not alter the exercise descriptions in their homework file, the student solutions are extracted and ran inside a fresh copy of the original Waterproof document.

\subsection{Evaluation}

We have evaluated Waterproof based on small student surveys, one-on-one conversations
with students, and conversations with the tutors who were able to help students with Waterproof exercises. Note that this is an evaluation of a previous version of Waterproof which was used in last year (2022-2023)'s Analysis 1 course. We expect that the updated version of the software presented in this article will have improved the user experience.

\subsubsection{Student Experience}

Student opinions on Waterproof vary: some prefer to do their homework by hand, some are set on completing all exercises in Waterproof; both groups seem to be evenly represented. 
Both stronger and weaker students use the software.
The retention rate is high: at the start of last year's course (2022-2023), 25 student groups ($\approx$\,100 students) handed in their homework using Waterproof; 19 groups ($\approx$\,76 students) continued using Waterproof for the final homework exercise. These approximate student numbers should be taken with a grain of salt, we did not track the number of individual students that used Waterproof. In some groups, a single group member works on all the Waterproof exercises, becoming this group's `Waterproof expert'; in other groups, all four students work on the Waterproof exercises together.

Students are sometimes frustrated with Waterproof, relating to aspects of the program that still need to be improved. Syntax errors, for example, provide too little information on how to fix them: a single, hard-to-notice typo can cause an entire proof step to be rejected. Sometimes also correct proof steps are rejected by Waterproof: some statements are too difficult for the automated proving system to verify, whereas they are `basic' for both students and instructors. Even in these cases, some students are determined to finish the exercise using Waterproof, even if it means rewriting their proof multiple times in different ways until it is accepted.

Students tend to stop using Waterproof when they feel they have little to gain from using the software, meaning that they think their handwritten solutions will get them enough points that working out their proofs in full using Waterproof is not worth the extra time and effort. Considering the high retention rate, this implies most students still find benefit to using Waterproof even at the end of the course.

\subsubsection{Teacher Experience}

The instructors who used Waterproof were impressed by how easy it was to use without prior experience with proof assistants. Nevertheless, they ran into the same issues as students, namely that some obvious statements would not be approved by the program and that some typos were hard to fix due to incomplete error messages.

Teachers also noticed that a majority of the questions they received during instruction classes were related to Waterproof instead of the exercises which were not converted into Waterproof documents. Students confirmed that they preferred to work on the Waterproof exercises during the instructions because this allowed them to ask questions directly when they encountered errors.

\subsubsection{Observations regarding Waterproof's Effect on Handwritten Proofs}

We observed that students started using Waterproof's  specific formulations of proof steps in their handwritten homework as well as on the exam. The explicit use of phrases like \emph{``We need to show that ...''}, \emph{``By~...~it holds that ...''}, or \emph{``We claim that ...''} clearly communicated students' intentions and made it easier to grade their proofs, even if they turned out to be incorrect. Variables and hypotheses were properly introduced before being used, and quantifiers were not kept around unnecessarily.

\section{Discussion}
Preliminary observations suggest that using Waterproof helps to improve the quality of students' handwritten proofs, but this claim needs to be backed up by a proper study. Although students were separated into instruction classes that offered Waterproof support and ones that did not, we did not take the preparatory steps, e.g. getting students' consent, needed for a publishable evaluation. These steps have been taken for this year (2023-2024)'s Analysis 1 course; we will study student performance based on their grades and the structure of their proofs.

Although Waterproof's custom proof language can be used across mathematical disciplines, creating exercises for different courses remains a time-consuming process for teachers and/or developers. Part of the problem is finding a suitable library wherein all the required mathematics have been formalized. Often these libraries encode definitions and theorems in different ways than the precise formulations used in a course, so a bridge needs to be built to connect the two. Students will get confused if definitions in Waterproof exercises differ from those used in the regular exercises. 
Besides finding/building the right mathematics library, using Waterproof's custom proof language comes with an additional requirement: having to tune its automated proof finding system. The right lemmas have to be provided to make sure that statements which are `obvious' to students and instructors are not rejected by the program.

The formulation of the proof steps in Waterproof's custom proof language are quite strict. Even small variations, like being able to write \colorbox{gray!20}{\texttt{Assume (...).}} instead of \colorbox{gray!20}{\texttt{Assume that (...).}}, have to be hard-coded. Deviations from these pre-defined formulations can result in syntax errors that are difficult to understand, this is especially a problem with typos. Currently, only English language proof step formulations are supported, other languages could be added without too much effort by developers.

Our solution for reasoning with chains of (in)equalities works very well. Still, considering how ubiquitous this reasoning technique is in mathematics, it would be nice if it were officially supported by \textsc{Coq}.

\subsection{Future Work}
We have iterated on Waterproof for the last four years and we plan to continue improving the software.

The automation system used by the custom proof language can still be improved further.
Last year, some statements which students and teachers considered obvious were still rejected by the program. We have further tuned the system's hint databases to prevent such cases from happening.
This year's iteration of the Analysis 1 course (2023--2024) will be another benchmark to see how successful we have been.

The feedback provided by the automation system when using the \colorbox{gray!20}{\texttt{By ...}}-clause can also be improved further. Currently, we are able to check whether the justification provided by the user is indeed necessary to proof their claim for external lemmas. It is, however, also possible to insert the label of local assumptions into the \colorbox{gray!20}{\texttt{By ...}}-clause, and these justifications cannot yet be checked. Moreover, if the specified lemma cannot be used, because, for example, one of its preconditions was not met, the user is only informed that the lemma failed to be used, and not about the reason why.

We also need better control over the automation system's ability to expand definitions. Currently, most of the tactics and solvers used by the automation system are able to do this automatically. As a result, some exercises can be solved almost directly by the automation system. This teaches wrong behavior to students since a proof might be accepted as complete without including some important intermediate steps. Ideally, the automatic expansion of some definitions should be shielded, such that the user has to explicitly mention the definition as part of a statement's justification, e.g. by writing
\begin{coqcode}
    By definition_of_supremum it suffices to show that (...).
\end{coqcode}
\vspace{-5pt}
    
We have not yet found a satisfying way to encode subsets in proof assistants that allows them to be used like subsets in regular mathematics. Both the implementation as sigma (or record) types and as classifying predicates have their disadvantages that we do not want to bother students with: sigma types require coercions for elements from the subset to be used as elements from the underlying set, and quantification over the elements satisfying a classifying predicate has to be done via quantification over the underlying set. The mathematical libraries we considered also did not stick to one approach, they mixed them. For example, although most libraries tend towards using classifying predicates, the subset relation $\mathbb N \subset \mathbb R$ is almost always implemented using a coercion.

Currently the editor only handles \texttt{.mv} \textsc{Coq} files; we plan to support the older \texttt{.v}-format in the future, rendering the \textsc{Coqdoc} comments as formatted text. We are also working to improve the modularity of the Waterproof editor, so others can reuse the parts that they like in their own projects. In the long term, we would like to have a version of Waterproof that runs completely in the browser like \textsc{jsCoq} \cite{jsCoq}.

We have chosen to, by default, only have the Waterproof editor show a restricted version of the automated bookkeeping overview, but others might disagree with this choice. Although a restricted overview helps with our goal to make the process of writing a proof in Waterproof more similar to writing a proof with pen-and-paper, displaying all the available information also has its educational benefits. For example, showing all the variables, assumptions, and intermediate assertions allows students to focus on how to connect the logical arguments, or it might help students to better understand what the effects of certain proof steps are. It would be interesting to explore the impact of limiting the scaffolding provided by the bookkeeping overview on the students' ability to transition from writing proofs in a proof assistant to writing proofs with pen-and-paper.

\section{Conclusion}

Waterproof is an adaptation of the \textsc{Coq} proof assistant designed for teaching how to write mathematical proofs. Although regular proof assistants, with their direct feedback capabilities, have the potential to serve as a training areas for proof writing, they are currently lacking for use in education: there is a steep learning curve and proof writing skills obtained in these programs do not seem to transfer to handwritten proofs. We attempt to address both issues by making writing a proof in Waterproof more closely resemble writing a proof by hand.

Waterproof consists of two components, the first of which is a custom proof language that allows \textsc{Coq} proofs to be written in a style more similar to handwritten proofs. Proof steps take the form of sentences used in regular mathematics and an implicit automation system is used to prevent users from having to provide justifications for basic statements. At crucial steps, users are forced to signpost their proofs to keep them readable. We also add the option to write chains of (in)equalities, a feature missing in default \textsc{Coq}.

The second component is a custom editor, which uses mixed documents that contain both verified \textsc{Coq} proofs and formatted text. Some of its features are designed specifically for use in education, like the ability to limit student input to certain parts of the document, and the by-default restriction of the information shown by the automated bookkeeping overview and the scaffolding this provides. The editor also comes with some quality-of-life features, like autocompletion for proof steps and mathematical symbols.

The Waterproof software has been successfully used to supplement the teaching of mathematical proof writing in the Analysis 1 course at the Eindhoven University of Technology.
Preliminary observations suggest that using Waterproof assists students in clarifying the logical structure of their proofs, including those written by hand.
A future study will investigate these claims more thoroughly by systematically analyzing student grades and the quality of their proofs.

\paragraph{Acknowledgements.}

In 2019, the first version of Waterproof arose as the result of a Software Engineering Project (SEP) at Eindhoven University of Technology, which is part of the bachelor curriculum. Iterations by other SEP Teams followed in 2021 and 2023. We would like to thank the members of the SEP teams ChefCoq (2019), Waterfowl (2021), KroQED (2023) and their supervisors: Thijs Beurskens, Bas Gieling, Stijn Gunter, Menno Hofst\'{e}, Hugo Melchers, Anne Nijsten, Reinier Schmiermann, Erik Takke, David Tuin, Ruben Verhaegh, Geert van Wordragen, and their supervisor Tom Verhoeff; Adrien Castella, Adrian Cuco\c{s}, Cosmin Manea, Noah van der Meer, Lulof Pir\'{e}e, Mihail \c{T}ifrea, Tristan Trouwen, Tudor Voicu, Adrian Ș.~Vrămuleț, Yuqing Zeng, and their supervisor Gerard Zwaan; Dick Arends, Franciska Asma, Casper Bloemhof, Daniel Chou Rainho, Milo Goolkate, Collin Harcarik, Banda Norimarna, Gijs Pennings, Lisa Verhoeven, Pieter Wils and their supervisor Roel Bloo. 
We thank Georgios Skantzaris for his feedback and help with testing when getting Waterproof ready for teaching for the first time.
We would like to thank EdIn (Education Innovation) at Eindhoven University of Technology for financing the work of the students Thijs Beurskens, Sean McCarren, Jan Moraal and David Tuin on this project.

We would like to thank Balthazar Patiachvili for improving the automation system and converting part of the coq-waterproof library to an OCaml plugin during his 2023 internship.
Finally, we are very grateful to Emilio Jesús Gallego Arias for his support, especially with \textsc{Coq SerAPI} and \textsc{Coq LSP}, and his quick answers to many questions.

\nocite{*}
\bibliographystyle{eptcs}
\bibliography{paper}
\end{document}